\documentclass[11pt]{amsart}
\usepackage{amssymb, mathtools, fullpage, color, ytableau, bbm}
\ytableausetup{smalltableaux}
\usepackage[colorinlistoftodos]{todonotes}
\usepackage[colorlinks=true, pdfstartview=FitV,linkcolor=blue,citecolor=blue,urlcolor=blue]{hyperref}

\theoremstyle{definition}
\newtheorem {theorem}{Theorem}[section]
\newtheorem {lemma}[theorem]{Lemma}

\newtheorem {conjecture}[theorem]{Conjecture}
\newtheorem{definition}[theorem]{Definition}
\newtheorem{remark}[theorem]{Remark}
\newtheorem{example}[theorem]{Example}

\newenvironment{red}{\relax\color{red}}{\relax}
\newenvironment{blue}{\relax\color{blue}}{\hspace*{.5ex}\relax}

\newcommand{\ber}{\begin{red}}
\newcommand{\er}{\end{red}}
\newcommand{\beb}{\begin{blue}}
\newcommand{\eb}{\end{blue}}

\newcommand{\seteq}{\coloneqq}

\newcommand{\GL}{\operatorname{GL}}

\numberwithin{equation}{section}
\renewcommand{\baselinestretch}{1.3}

\begin{document}

\title[Data-scientific study of  Kronecker coefficients]{Data-scientific study of  Kronecker coefficients}

\date{\today}

\author[K.-H. Lee]{Kyu-Hwan Lee}
\address{Department of Mathematics, University of Connecticut, Storrs, CT 06269, U.S.A.}
\email{khlee@math.uconn.edu}

\keywords{Kronecker coefficients, principal component analysis, kernel methods}

\begin{abstract}
We take a data-scientific approach to study whether Kronecker coefficients are zero or not. Motivated by principal component analysis and kernel methods, we define {\em loadings} of partitions and use them to describe a sufficient condition for Kronecker coefficients to be nonzero. The results provide new methods and perspectives for the study of these coefficients.  
 
\end{abstract}

\maketitle

\section{Introduction}\label{sec1}

In recent years, there has been much discussion about the potential for AI and machine learning to change mathematics research (e.g. \cite{DVB+,W,B}). Numerous examples demonstrate machine learning's ability to discern patterns in mathematical datasets (e.g. \cite{Ashmore:2019wzb,Carifio:2017bov,DLQ,He:2019nzx,HLOa,HLOb,HLOc,Jejjala:2019kio}). The recent discovery \cite{HLOP} of a new phenomenon, called {\em murmuration},  illustrates the significant potential of considering mathematical objects within the framework of data science. All these developments remind us of how fruitful it has been to study datasets in modern mathematics. For instance, the prime number theorem and the Birch--Swinnerton-Dyer conjecture are results of investigating certain datasets.  

The goal of this paper is to apply this paradigm of mathematics research to representation theory. One of the primary objectives in representation theory is to decompose a representation into its irreducible components, with algebraic combinatorics providing a vital and practical method for describing this decomposition. A prototypical example is the decomposition of the tensor product of two irreducible representations of the general linear group $\GL_N(\mathbb C)$, where the Littlewood--Richardson rule completely describes the decomposition using skew semi-standard tableaux.
 
Surprisingly enough, there has been no similar success with the symmetric group until now.  
Let $\mathfrak S_n$ be the symmetric group of degree $n$. The irreducible representations $S_\lambda$ of $\mathfrak S_n$ over $\mathbb C$ are parametrized by partitions $\lambda$ of $n$, written as $\lambda \vdash n$. The tensor product of two irreducible representations $S_\lambda$ and $S_\mu$ ($\lambda, \mu \vdash n$) is decomposed into a sum of irreducible representations:
\[  S_\lambda \otimes S_\mu = \bigoplus_{\nu\vdash n} g_{\lambda, \mu}^\nu S_\nu  \quad (g_{\lambda, \mu}^\nu \in \mathbb Z_{\ge 0}) . \]
The decomposition multiplicities $g_{\lambda, \mu}^\nu$ are called the {\em Kronecker coefficients}.

In stark contrast to the Littlewood--Richardson coefficients for $\GL_N(\mathbb C)$, no combinatorial description has been known for $g_{\lambda, \mu}^\nu$ since Murnaghan \cite{Mur} initially posed the question in 1938, and it is still considered as one of the main problems in the combinatorial representation theory. Only partial results are available due to Remmel \cite{Rem}, Ballantine--Orellana \cite{BO}, Remmel--Whitehead \cite{RW}, Blasiak--Mulmuley--Sohoni \cite{BMS} and Blasiak \cite{Bla}. Recently, the coefficients $g_{\lambda, \mu}^\nu$ have also been studied from the perspective of computational complexity. Notably,  Ikenmeyer, Mulmuley and Walter \cite{IMW} demonstrated that determining  whether a given Kronecker coefficient is non-zero is NP-hard. Additionally, Pak and Panova \cite{PP1, PP2} have made other significant contributions to this topic.

In the previous article \cite{L}, we applied standard machine learning tools to datasets of the Kronecker coefficients, and observed that the trained classifiers attained high accuracies ($> 98\%$) in determining whether Kronecker coefficients are zero or not. The outcomes clearly suggest that further data-scientific analysis may reveal new structures in the datasets of the Kronecker coefficients. 
In this paper, we indeed find new structures; more precisely, we adopt ideas from principal component analysis (PCA) and kernel methods to define the {\em similitude} matrix and the {\em difference} matrix for the set $\mathcal P(n)$ of partitions of $n$.  Then we introduce {\em loadings} of the partitions in terms of eigenvectors associated to the largest eigenvalues of these matrices, and use the loadings to describe a sufficient condition for the Kronecker coefficients to be nonzero. This condition can be used very effectively. See \eqref{eqn-box} and Example \ref{exa-18} below it.

The observations made in this paper are purely data-scientific and experimental, and no attempts are undertaken to prove them using representation theory. Also, it should be noted that our sufficient condition does not cover the {\em middle part} where loadings for zero and nonzero Kronecker coefficients overlap. Since our method is a variation of PCA, it is essentially linear. In order to cover the middle part, it is likely that one needs to take   a nonlinear approach. The aforementioned high accuracies reported in \cite{L} indicate that efficient strategies may be developed to go much deeper into the middle part.    

After this introduction, in Section \ref{s:Ltraining}, we define the similitude and difference matrices and the loadings of partitions. In Section \ref{sec-loadings}, we investigate the probabilistic distributions of loadings. In the final section, we consider the minimum values of the loadings to determine whether the Kronecker coefficients are zero or nonzero. In Appendix, we tabulate the loadings of partitions in $\mathcal P(n)$ for $6 \le n \le 12$.

\subsection*{Acknowledgments}
 The author is grateful to Alex Davies and Been Kim for their helpful discussions. He also thanks the anonymous referee whose comments significantly improved the exposition of this paper.  The author would like to acknowledge the support and hospitality of the Isaac Newton Institute for Mathematical Sciences and the Center for Quantum Structures in Modules and Spaces. This work was partially supported by an EPSRC grant \#EP/R014604/1 and by a grant from the Simons Foundation (\#712100).

\section{Similitude and difference matrices}\label{s:Ltraining}

Let $\mathfrak S_n$ be the symmetric group of degree $n$ and consider representations of $\mathfrak S_n$ over $\mathbb C$. The irreducible representations $S_\lambda$ of $\mathfrak S_n$ are parametrized by partitions $\lambda \in \mathcal P(n)$. Consider the tensor product of two irreducible representations $S_\lambda$ and $S_\mu$ for $\lambda, \mu \in \mathcal P(n)$. Then the tensor product is decomposed into a sum of irreducible representations:
\[  S_\lambda \otimes S_\mu = \bigoplus_{\nu\vdash n} g_{\lambda, \mu}^\nu S_\nu  \quad (g_{\lambda, \mu}^\nu \in \mathbb Z_{\ge 0}) . \]
The decomposition multiplicities $g_{\lambda, \mu}^\nu$ are called the {\em Kronecker coefficients}.

\medskip

There are symmetries among $g_{\lambda, \mu}^\nu$.
\begin{lemma} \cite[p.61]{FH} \label{lem-perm}
Let $\lambda, \mu, \nu \vdash n$. Then the Kronecker coefficients $g_{\lambda, \mu}^\nu$ are invariant under the permutations of $\lambda, \mu, \nu$. That is, we have
\[ g_{\lambda,  \mu}^\nu=g_{\mu, \lambda}^\nu=g_{\lambda, \nu}^\mu=g_{\nu, \lambda}^\mu=g_{\mu, \nu}^\lambda=g_{\nu, \mu}^\lambda. \]
\end{lemma}
For a partition $\lambda=(\lambda_1 \ge \lambda_2 \ge \cdots) $ of $n$, define $d_\lambda \seteq n- \lambda_1$, called the {\em depth} of $\lambda$. The following theorem provides a necessary condition for the Kronecker coefficient $g_{\lambda,\mu}^\nu$ to be nonzero. Other necessary conditions for $g_{\lambda,\mu}^\nu \neq 0$, which generalize Horn inequalities, can be found in \cite{Res}. We will describe a sufficient condition for for $g_{\lambda, \mu}^\nu \neq 0$ in this paper.  

\begin{theorem} \cite[Theorem 2.9.22]{JK} \label{thm-lim}
If $g_{\lambda, \mu}^\nu \neq 0$ then  \begin{equation} \label{eqn-dl}  |d_\lambda - d_\mu | \le d_\nu \le d_\lambda + d_\mu . \end{equation}
\end{theorem}

\bigskip

Now, for $n \in \mathbb Z_{>0}$, let $\mathcal P(n)$ be the set of partitions of $n$ as before.
We identify each element $\lambda$ of $\mathcal P(n)$ with a sequence of length $n$  by appending as many $0$-entries as needed.     We also consider $\mathcal P(n)$ as an ordered set by the lexicographic order. 

\begin{example} When $n=6$, we have
\begin{align*} \mathcal P(6) = \{ (6,0,0,0,0,0), (5, 1,0,0,0,0), (4, 2,0,0,0,0), (4, 1, 1,0,0,0),\\ (3,3,0,0,0,0), (3, 2, 1,0,0,0), (3, 1, 1, 1,0,0),(2,2,2,0,0,0),\\ (2, 2, 1, 1,0,0), (2, 1, 1, 1,1,0), (1,1, 1, 1, 1, 1) \} . \phantom{LLLLLLLLa} \end{align*}
\end{example}

When there is no peril of confusion, we will skip writing 0's in the sequence  when it ends with a string of 0's.  For instance, we write $(5,1)$ for $(5,1,0,0,0,0)$. Moreover, when the same part is repeated multiple times we may abbreviate it into an exponent. For example, $(2,1,1,1,1,1)$ may be written as $(2,1^5)$. 
The size of the set $\mathcal P(n)$ will be denoted by $p(n)$, and the set of triples $\mathbf t = (\lambda, \mu, \nu)$ of partitions of $n$ will be  denoted  by $\mathcal P(n)^3 \seteq \mathcal P(n) \times \mathcal P(n) \times \mathcal P(n)$. A partition is depicted by a collection of left-justified rows of boxes. For example, partition $(5,4,1)$ is depicted by $\raisebox{6 pt}{\ydiagram{5,4,1}}$.
The {\em conjugate} or {\em transpose} of a partition is defined to be the flip of the original diagram along the main diagonal. Hence the conjugate of $(5,4,1)$ is $(3,2,2,2,1)$ as you can see below:
\[  \ydiagram{5,4,1} \qquad \longleftrightarrow \qquad \raisebox{10 pt}{\ydiagram{3,2,2,2,1}} \]
 
\bigskip

Let $\mathsf P_n$ be the $p(n) \times n$ matrix having elements of $\mathcal P(n)$ as rows, and define the $p(n) \times p(n)$ symmetric matrix \[ \mathsf Y_n \seteq \mathsf P_n \mathsf P_n^{\top}.\] The matrix $\mathsf Y_n$ will be called the {\em similitude} matrix of $\mathcal P(n)$.

\begin{example} When $n=6$, we obtain 
\[{\scriptsize \mathsf P_6= \begin{bmatrix} 6&0&0&0&0&0 \\ 5& 1&0&0&0&0\\4& 2&0&0&0&0\\ 4& 1& 1&0&0&0\\ 3&3&0&0&0&0\\ 3& 2& 1&0&0&0\\3& 1& 1& 1&0&0\\2&2&2&0&0&0\\ 2& 2& 1& 1&0&0\\2& 1& 1& 1&1&0\\1&1& 1& 1& 1& 1\end{bmatrix} } \quad \text{ and } 
\quad {\scriptsize \mathsf Y_6=\left [ \begin{array}{ccccccccccc}36&30&24&24&18&18&18&12&12&12& 6 \\ 30& 26& 22& 21& 18& 17& 16& 12& 12& 11& 6\\ 24& 22& 20& 18& 18& 16& 14& 12& 12& 10& 6\\ 24& 21& 18& 18& 15& 15& 14& 12& 11& 10& 6\\ 18& 18& 18& 15& 18& 15& 12& 12& 12&  9& 6\\ 18& 17& 16& 15& 15& 14& 12& 12& 11&  9& 6\\ 18& 16& 14& 14& 12& 12& 12& 10& 10&  9& 6\\ 12& 12& 12& 12& 12& 12& 10& 12& 10&  8& 6\\ 12& 12& 12& 11& 12& 11& 10& 10& 10&  8& 6\\ 12& 11& 10& 10&  9&  9&  9&  8&  8&  8& 6\\ 6&  6&  6&  6&  6&  6&  6&  6&  6&  6& 6 \end{array} \right ] }.\]
Note that an entry $y_{\lambda, \mu}$ of $\mathsf Y_n=[y_{\lambda, \mu}]$ is indexed by $\lambda, \mu \in \mathcal P(n)$.
\end{example}

Since the matrix $\mathsf Y_n$ is symmetric, all its eigenvalues are real. Moreover, the Perron--Frobenius theorem \cite[Section III.2]{Ga} tells us that $\mathsf Y_n$ has a unique eigenvalue of largest magnitude and that the corresponding eigenvector can be chosen to have strictly positive components.

\begin{definition} \label{def-p}
Let $\mathbf v =(v_\lambda)_{\lambda \in \mathcal P(n)}$ be an eigenvector of the largest eigenvalue of $\mathsf Y_n$ such that $v_\lambda >0$ for all $\lambda \in \mathcal P(n)$. Denote by $v_{\mathrm{max}}$ (resp. $v_{\mathrm{min}}$) a maximum (resp. minimum) of $\{ v_\lambda \}_{\lambda \in \mathcal P(n)}$. Define \[ a_\lambda \seteq 100 \times \frac{v_\lambda - v_{\mathrm{min}}}{v_{\mathrm{max}} - v_{\mathrm{min}}} \quad \text{for } \lambda \in \mathcal P(n) .\] The value $a_\lambda$ is called the {\em $a$-loading} of partition $\lambda \in \mathcal P(n)$. 
\end{definition}

The above definition presents a novel concept in the exploration of Kronecker coefficients. However, when $n$ is large, one might wonder how to compute an eigenvector of the largest eigenvalue of $\mathsf Y_n$. The direct computation of eigenvalues can be computationally intensive. For an $N \times N$ matrix, it is known that direct computations of eigenvalues have a time complexity of $O(N^3)$.

An efficient algorithm to calculate an eigenvector $\mathbf v$ in Definition \ref{def-p} is the {\em power iteration}: Let $\mathbf v_0=(1,0, \dots , 0)^\top$ be the first standard column vector. Inductively, for $k=0,1,2,\dots$, define
\[ \mathbf v_{k+1} = \frac{\mathsf Y_n \mathbf v_k}{ \lVert \mathsf Y_n \mathbf v_k \rVert_2} , \] where $\lVert (x_1, x_2, \dots, x_n)^\top \rVert_2 = ( \sum_{i=1}^n x_i^2)^{1/2}$. Then the limit
\[ \mathbf v = \lim_{k \to \infty} \mathbf v_k \] is an eigenvector of the largest eigenvalue of $\mathsf Y_n$.

\begin{example}
When $n=6$,  we have
{\scriptsize \begin{align*}
\mathbf v_1 & \approx (0.5203, 0.4336 ,0.3468 ,0.3468 ,0.2601 ,0.2601 ,0.2601 ,0.1734 ,0.1734 ,0.1734 ,0.0867)^\top , \\
\mathbf v_2 & \approx (0.4514 ,0.4022 ,0.3530  ,0.3377 ,0.3038 ,0.2885 ,0.2670 ,0.2240 ,0.2178 ,0.1934 ,0.1188)^\top , \\
\mathbf v_3 & \approx (0.4441 ,0.3985 ,0.3530  ,0.3366 ,0.3074 ,0.2910 ,0.2678 ,0.2291 ,0.2222 ,0.1957 ,0.1225)^\top , \\
\mathbf v_4 & \approx (0.4434 ,0.3982 ,0.3529 ,0.3365 ,0.3077,0.2913 ,0.2678 ,0.2296 ,0.2226 ,0.1960,0.1229)^\top , \\ 
\mathbf v_5 & \approx (0.4433 ,0.3981 ,0.3529 ,0.3365 ,0.3077,0.2913 ,0.2678 ,0.2297 ,0.2226 ,0.1960,0.1229)^\top , \\
\mathbf v_6 & \approx (0.4433 ,0.3981 ,0.3529 ,0.3365 ,0.3077,0.2913 ,0.2678 ,0.2297 ,0.2227 ,0.1960,0.1229)^\top.
\end{align*} }
Thus we can take as an approximation {\scriptsize \[ \mathbf v= (0.4433,0.3981,0.3529,0.3365,0.3077,0.2913,0.2678,0.2297,0.2227,0.1960,0.1229)^\top , \] }  and  the $a$-loadings are given by
\[ (a_\lambda)_{\lambda \in \mathcal P(n)} = (100.00,85.89,71.79,66.66,57.68,52.55,45.23,33.32,31.12,22.81,0.00).\] 
\end{example}

In the above example of $n=6$, we see that the $a$-loadings are compatible with the lexicographic order. In particular, the partition $(6)$ has $a$-loading $100$ and $(1,1,1,1,1,1)$ has $a$-loading $0$.
However, in general, the $a$-loadings are {\em not completely} compatible with the lexicographic order though they are strongly correlated. For instance, when $n=9$, the partition $(5,1,1,1,1)$ has $a$-loading $55.32$, while $(4,4,1)$ has $56.55$. See Appendix \ref{appen} for the values of $a$-loadings.  On the other hand, we say that $\lambda = (\lambda_1 \ge \lambda_2 \ge \dots \ge \lambda_n)$ dominates $\mu= (\mu_1 \ge \mu_2 \ge \dots \ge \mu_n)$ in the dominance order if $\lambda_1 + \cdots + \lambda_k \ge \mu_1 + \cdots \mu_k$ for all $k \ge 1$. Now one can observe that {\em the $a$-loadings are compatible with the dominance order.} \footnote{This was noticed by David Anderson after the first version of this paper was posted on the arXiv.} 

\bigskip

Define a $p(n) \times p(n)$ symmetric matrix $\mathsf Z_n=[z_{\lambda, \mu}]_{\lambda, \mu \in \mathcal P(n)}$ by
\[ z_{\lambda, \mu}= \lVert \lambda-\mu \rVert_1 \seteq \sum_{i=1}^n | \lambda_i - \mu_i| \]
for $\lambda=(\lambda_1, \lambda_2, \dots , \lambda_n)$ and $\mu=(\mu_1,\mu_2, \dots , \mu_n) \in \mathcal P(n)$.
The matrix $\mathsf Z_n$ will be called the {\em difference} matrix of $\mathcal P(n)$.

\begin{example}
When $n=6$, we obtain 
\[{\scriptsize \mathsf Z_6 = \left [ \begin{array}{ccccccccccc}  0  &2  &4  &4  &6  &6  &6  &8  &8  &8 &10 \\ 2  &0  &2  &2  &4  &4  &4  &6  &6  &6  &8 \\  4  &2  &0  &2  &2  &2  &4  &4  &4  &6  &8\\4  &2  &2  &0  &4  &2  &2  &4  &4  &4  &6\\6  &4  &2  &4  &0  &2  &4  &4  &4  &6  &8\\6  &4  &2  &2  &2  &0  &2  &2  &2  &4  &6\\6  &4  &4  &2  &4  &2  &0  &4  &2  &2  &4\\8  &6  &4  &4  &4  &2  &4  &0  &2  &4  &6\\8  &6  &4  &4  &4  &2  &2  &2  &0  &2  &4\\8  &6  &6  &4  &6  &4  &2  &4  &2  &0  &2\\ 10  &8  &8  &6  &8  &6  &4  &6  &4  &2  &0 \end{array} \right ] }.\]  
\end{example}

Similarly to $\mathsf Y_n$, all the eigenvalues of $\mathsf Z_n$ are real. It is easy to see that $\mathsf Z_n$ is irreducible, and hence the Perron--Frobenius theorem for matrices with nonnegative entries \cite[Section III.2]{Ga} tells us that $\mathsf Z_n$ has a unique eigenvalue of largest magnitude and that the corresponding eigenvector can be chosen to have strictly positive components.

\begin{definition}
Let $\mathbf w =(w_\lambda)_{\lambda \in \mathcal P(n)}$ be an eigenvector of the largest eigenvalue of $\mathsf Z_n$ such that $w_\lambda >0$ for all $\lambda \in \mathcal P(n)$. Denote by $w_{\mathrm{max}}$ (resp. $w_{\mathrm{min}}$) a maximum (resp. minimum) of $\{ w_\lambda \}_{\lambda \in \mathcal P(n)}$. Define \[ b_\lambda \seteq 100 \times \frac{w_\lambda - w_{\mathrm{min}}}{w_{\mathrm{max}} - w_{\mathrm{min}}} \quad \text{for } \lambda \in \mathcal P(n) .\] The value $b_\lambda$ is called the {\em $b$-loading} of partition $\lambda \in \mathcal P(n)$. 
\end{definition}

Like Definition \ref{def-p}, the above definition introduces a new concept into the study of Kronecker coefficients.  We will show its usefulness in Section \ref{sec-noe}. 
The power iteration works equally well to compute $\mathbf w$: Let $\mathbf w_0=(1,0, \dots , 0)^\top$ and define
\[ \mathbf w_{k+1} = \frac{\mathsf Z_n \mathbf w_k}{ \lVert \mathsf Z_n \mathbf w_k \rVert_2} . \] Then the limit
\[ \mathbf w = \lim_{k \to \infty} \mathbf w_k \] is an eigenvector of the largest eigenvalue of $\mathsf Z_n$.  

\medskip

\begin{example}
When $n=6$, we have
{\scriptsize \begin{align*}
\mathbf w_1& \approx (0.0000,     0.0958,0.1916,0.1916,0.2873,0.2873,0.2873,0.3831,0.3831,0.3831,0.4789)^\top , \\
\mathbf w_2 & \approx  (0.5177,0.3705,0.2992,0.2565,0.3087,0.2042,0.2042,0.2517,0.1947,0.2280,0.3277)^\top , \\
\vdots &  \\
\mathbf w_{10} & \approx  (0.4046,0.2962,0.2662,0.2394,0.3061,0.2318,0.2393,0.3060,  0.2662,0.2961,0.4044)^\top , \\
\mathbf w_{11} & \approx  (0.4045,0.2961,0.2662,0.2393,0.3061,0.2318,0.2393,0.3061,0.2662,0.2962,0.4045)^\top , \\
\mathbf w_{12} & \approx  (0.4045,0.2961,0.2662,0.2393,0.3061,0.2318,0.2393,0.3061,0.2662,0.2961, 0.4045)^\top .
\end{align*} }
Thus we can take as an approximation  {\scriptsize \[ \mathbf w= (0.4045,0.2961,0.2662,0.2393,0.3061,0.2318,0.2393,0.3061,0.2662,0.2961, 0.4045)^\top , \] }  and  the $b$-loadings are given by
\[ (b_\lambda)_{\lambda \in \mathcal P(n)} = (100.00,37.25,19.93,4.36,43.01,0.00,4.36,43.01,19.93,37.25, 100.00).\]
\end{example}

\begin{remark} \label{quote-1}
In the above example, we notice that the partitions $(6,0,0,0,0,0)$  and $(1,1,1,1,1,1)$ both have $b$-loading $100$ and the partition $(3,2,1,0,0,0)$ has $b$-loading $0$.  In general, we observe that {\em if $\lambda$ and $\mu$ are conjugate in $\mathcal P(n)$, then their $b$-loadings are the same, i.e., $ b_\lambda = b_\mu$.}   
\end{remark}

\begin{remark}
It would be interesting to combinatorially characterize the loadings of $\lambda \in \mathcal P(n)$. 
\end{remark}

\bigskip

For $\mathbf t=(\lambda, \mu, \nu) \in \mathcal P(n)^3$, we will write
\[ g(\mathbf t) \seteq g_{\lambda, \mu}^\nu . \]

\begin{definition} \label{def-three}
Let $\mathbf t=(\lambda, \mu, \nu) \in \mathcal P(n)^3$. Define the {\em $a$-loading} (resp.  {\em $b$-loading})  of $\mathbf t$, denoted by  $a(\mathbf t)$ (resp. $b(\mathbf t)$), to be the sum of the $a$-loadings (resp. $b$-loadings) of $\lambda, \mu$ and $\nu$, i.e., 
\[ a(\mathbf t) \seteq a_\lambda + a_\mu + a_\nu \qquad (\text{resp.} \ b(\mathbf t) \seteq b_\lambda + b_\mu + b_\nu ) . \] 
\end{definition}

\subsection{Connections to PCA and  kernel methods }
The definitions of similitude and difference matrices are closely related to PCA and kernel methods (see, e.g.,  \cite[Sections 3.5 and 12.3]{Hastie}), respectively. Indeed, we look at the matrix $\mathsf P_n^\top$ as a data matrix.

\begin{example}
When $n=6$, we get 
\[{\scriptsize \mathsf P_6^\top = \left [ \begin{array}{rrrrrrrrrrr}
6 & 5 & 4 & 4 & 3 & 3 & 3 & 2 & 2 & 2 & 1 \\
0 & 1 & 2 & 1 & 3 & 2 & 1 & 2 & 2 & 1 & 1 \\
0 & 0 & 0 & 1 & 0 & 1 & 1 & 2 & 1 & 1 & 1 \\
0 & 0 & 0 & 0 & 0 & 0 & 1 & 0 & 1 & 1 & 1 \\
0 & 0 & 0 & 0 & 0 & 0 & 0 & 0 & 0 & 1 & 1 \\
0 & 0 & 0 & 0 & 0 & 0 & 0 & 0 & 0 & 0 & 1
\end{array}\right ]}, \] 
and consider this as a data matrix of 6 data points with 11 features. 
\end{example}

Since the average of each column is $1$ for $\mathsf P_n^\top$, the covariance matrix of the data matrix $\mathsf P_n^\top$ is $(\mathsf P_n - \mathbbm{1}) (\mathsf P_n - \mathbbm{1})^\top$, where $\mathbbm{1}$ is the matrix with all entries equal to $1$.  As there are no significant differences in the largest eigenvalues or the directions of their eigenvectors, we use the similitude matrix 
$\mathsf Y_n=\mathsf P_n \mathsf P_n^\top$ 
as a substitute for the covariance matrix. 
Then an eigenvector of the largest eigenvalue of $\mathsf Y_n$ provides a good approximation to a weight vector of the first principal component, leading to the definition of $a$-loadings.   

The idea of a kernel method is to embed a dataset into a different space of (usually) higher dimension. In order to utilize this idea, we consider the matrix $\mathsf P_n$ as a data matrix with $p(n)$ data points and $n$ features. Then we map a partition $\lambda$, which is an $n$-dimensional row vector of $\mathsf P_n$, onto the $p(n)$-dimensional vector $(\lVert \lambda -\mu \rVert_1)_{\mu \in \mathcal P(n)}$, and the resulting new matrix is exactly the difference matrix $\mathsf Z_n$.

\begin{example}
When $n=6$, we obtain 
\[{\scriptsize \mathsf P_6= \begin{bmatrix} 6&0&0&0&0&0 \\ 5& 1&0&0&0&0\\4& 2&0&0&0&0\\ 4& 1& 1&0&0&0\\ 3&3&0&0&0&0\\ 3& 2& 1&0&0&0\\3& 1& 1& 1&0&0\\2&2&2&0&0&0\\ 2& 2& 1& 1&0&0\\2& 1& 1& 1&1&0\\1&1& 1& 1& 1& 1\end{bmatrix}  \quad  \mapsto \quad  \mathsf Z_6 = \left [ \begin{array}{ccccccccccc}  0  &2  &4  &4  &6  &6  &6  &8  &8  &8 &10 \\ 2  &0  &2  &2  &4  &4  &4  &6  &6  &6  &8 \\  4  &2  &0  &2  &2  &2  &4  &4  &4  &6  &8\\4  &2  &2  &0  &4  &2  &2  &4  &4  &4  &6\\6  &4  &2  &4  &0  &2  &4  &4  &4  &6  &8\\6  &4  &2  &2  &2  &0  &2  &2  &2  &4  &6\\6  &4  &4  &2  &4  &2  &0  &4  &2  &2  &4\\8  &6  &4  &4  &4  &2  &4  &0  &2  &4  &6\\8  &6  &4  &4  &4  &2  &2  &2  &0  &2  &4\\8  &6  &6  &4  &6  &4  &2  &4  &2  &0  &2\\ 10  &8  &8  &6  &8  &6  &4  &6  &4  &2  &0 \end{array} \right ] }.\]  
\end{example}

Since the difference matrix $\mathsf Z_n$ is a symmetric matrix, we consider an eigenvector of the largest eigenvalue of $\mathsf Z_n$ to obtain the direction of largest variations in the differences. This leads to the definition of $b$-loadings.

\section{Distributions of loadings} \label{sec-loadings}
In this section, we present the histograms of loadings and describe the corresponding distributions. First, we consider all the triples of $\mathbf t \in \mathcal P(n)^3$, and after that, separate them according to whether $g(\mathbf t) \neq 0$ or $= 0$. 

\medskip

\begin{figure}[h!]
\begin{center}
\includegraphics[scale=0.35]{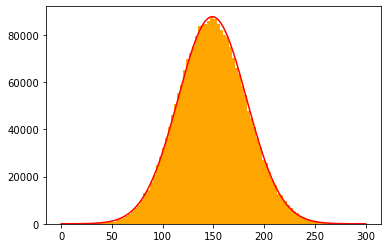}\quad \includegraphics[scale=0.35]{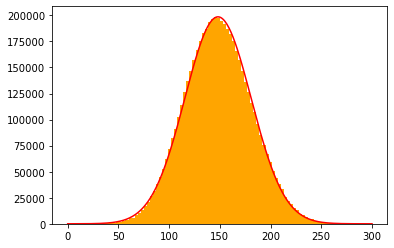}\quad \includegraphics[scale=0.35]{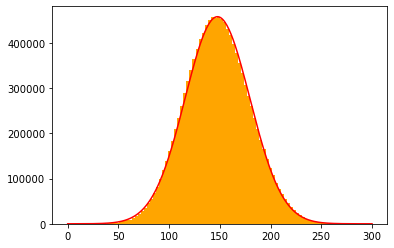}
\end{center}
\caption{\small\sf Histograms of $a$-loadings of $\mathbf t \in \mathcal P(n)^3$ for $n=14,15,16$ from left to right along with curves (red) of normal distributions} \label{fig-6} \end{figure}

\begin{figure}[h!]
\begin{center}
\includegraphics[scale=0.35]{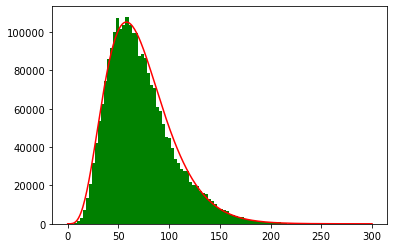}\quad \includegraphics[scale=0.35]{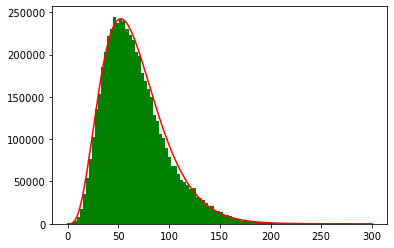}\quad \includegraphics[scale=0.35]{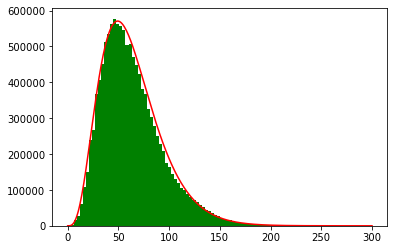}
\end{center}
\caption{\small\sf Histograms of $b$-loadings of $\mathbf t \in \mathcal P(n)^3$ for $n=14,15,16$ from left to right along with curves (red) of gamma distributions} \label{fig-7} \end{figure}

Figure \ref{fig-6} (resp. Figure \ref{fig-7}) has the histograms of $a$-loadings (resp. $b$-loadings) of $\mathbf t \in \mathcal P(n)^3$ for $n=14,15,16$. According to what the histograms suggest, we propose a conjecture: 
\begin{conjecture}
{\em Consider $\mathcal P(n)^3$ as a sample space. Then the sequence of random variables $X^a_n$ (resp. $X^b_n$) defined by the $a$-loadings (resp. $b$-loadings) of $\mathbf t$ converges in distribution to a normal (resp. gamma) random variable as $n \to \infty$}. 
\end{conjecture}
We sketch the curves of normal distributions on the histograms in Figure \ref{fig-6}. Here we note that the mean is not exactly 150. Actually, the mean values of the $a$-loadings are $\approx 148.86, 148.15, 147.65$ for $n=14,15,16$, respectively.
Similarly, we draw the curves of gamma distributions in Figure \ref{fig-7}. The mean values of the $b$-loadings are $\approx 72.07, 66.71, 63.48$ for $n=14,15,16$, respectively. 

\medskip

\begin{figure}[h!]
\begin{center}
\includegraphics[scale=0.35]{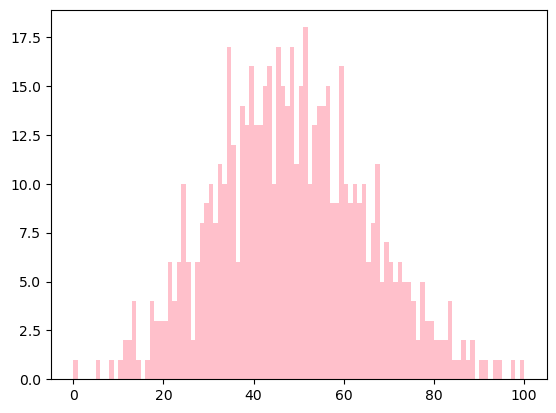}\quad \includegraphics[scale=0.35]{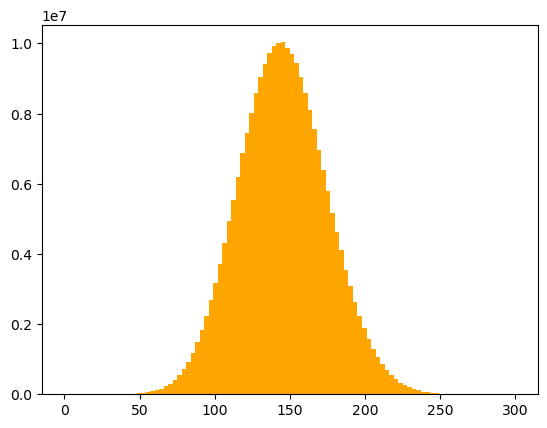}\\ \includegraphics[scale=0.35]{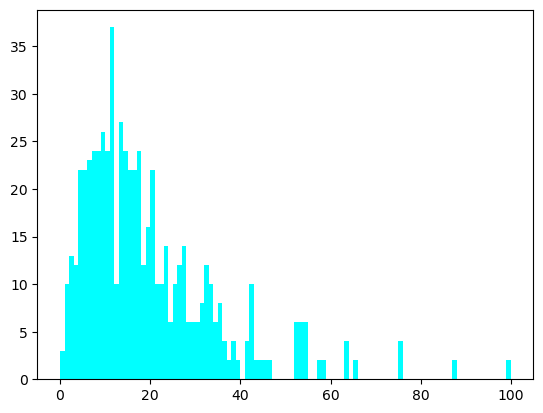} \quad \includegraphics[scale=0.35]{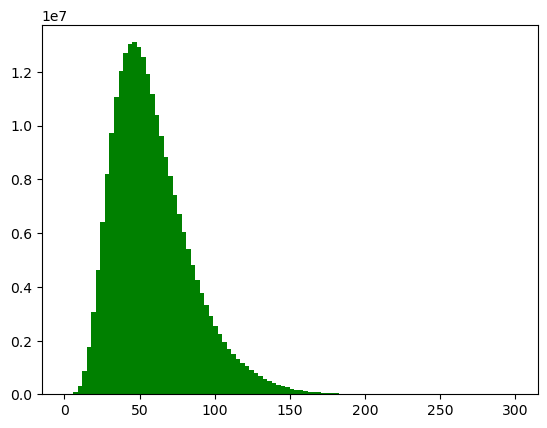}
\end{center}
\caption{\small\sf Histograms of $a$-loadings of $\lambda \in \mathcal P(n)$ (top-left) and $\mathbf t \in \mathcal P(n)^3$ (top-right) and histograms of $b$-loadings of $\lambda$ (bottom-left) and $\mathbf t$ (bottom-right) when $n=20$} \label{fig-8} \end{figure}

When $n=14,15,16$, the histograms of the loadings of partitions $\lambda \in \mathcal P(n)$,  in contrast to triples $\mathbf t \in \mathcal P(n)^3$, lack sufficient data points to ascertain their underlying distributions.   (Note that $p(16)= 231$.) Nonetheless,  since $a_\lambda$, $a_\mu$ and $a_\nu$ are computed independently for $a(\mathbf t) = a_\lambda +a_\mu + a_\nu$, and the sum of normal random variables is itself normal, it seems reasonable to expect that the $a$-loadings of $\lambda$ follow a normal distribution. By the same reasoning, we conjecture that the $b$-loadings of $\lambda$ follow a gamma distribution. If the conjectures are true, the loadings of $\mathbf t \in \mathcal P(n)^3$ will naturally have the distribution given as a sum of three independent distributions.  (Recall Definition \ref{def-three}.)   
Figure \ref{fig-8} has the histograms of loadings of $\lambda$ and $\mathbf t$ when $n=20$, which seem to be consistent with this expectation.


\section{Separation of $g(\mathbf t) \neq 0$ from $g(\mathbf t) = 0$} \label{sec-noe} 
In this section, we consider the distributions of loadings according to whether the Kronecker coefficients $g(\mathbf t)$ are zero or nonzero. Using minimum values of loadings in each case, we will obtain vertical lines which separate the distributions of these two cases. 

\medskip

\begin{figure}[h!]
\begin{center}
\includegraphics[scale=0.5]{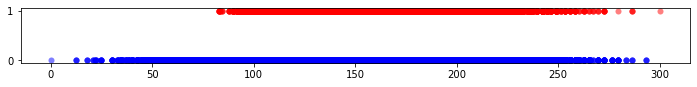}\\ \includegraphics[scale=0.5]{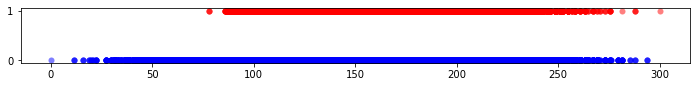}\\\includegraphics[scale=0.5]{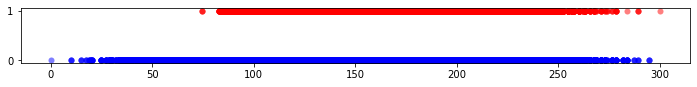}
\\\includegraphics[scale=0.5]{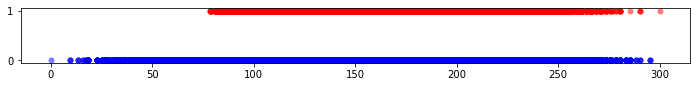}
\end{center}
\caption{\small\sf Ranges of $a$-loadings for $n=10,11,12,13$ from top to bottom. A red (resp. blue) dot at $(x,1)$ (resp. $(x,0)$) corresponds to $\mathbf t \in \mathcal P(n)^3$ with $a(\mathbf t) =x$ and $g(\mathbf t) \neq 0$ (resp. $g(\mathbf t)=0$).} \label{fig-4} \end{figure}

\begin{figure}[h!]
\begin{center}
\includegraphics[scale=0.35]{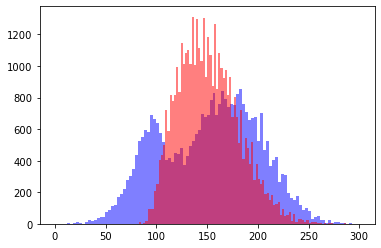}\quad \includegraphics[scale=0.35]{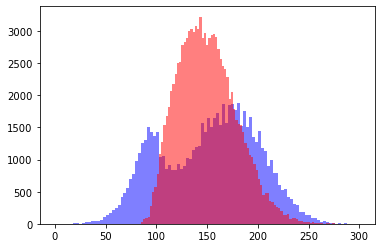}\\ \includegraphics[scale=0.35]{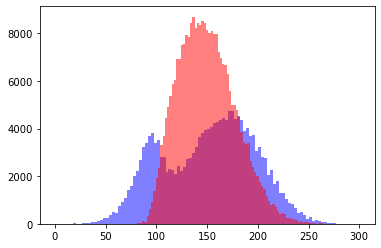} \quad \includegraphics[scale=0.35]{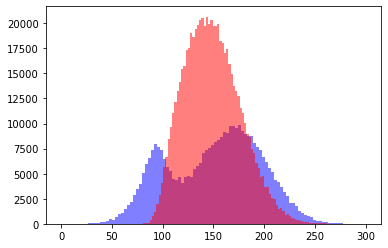}
\end{center}
\caption{\small \sf Histograms of $a$-loadings for $n=10$ (top-left), $11$ (top-right), $12$ (bottom-left) and $13$ (bottom-right). The red (resp. blue) region represents the numbers of $\mathbf t$ such that $g(\mathbf t) \neq 0$ (resp. $g(\mathbf t) =0$).} \label{fig-5} \end{figure}

\begin{figure}[h!]
\begin{center}
\includegraphics[scale=0.5]{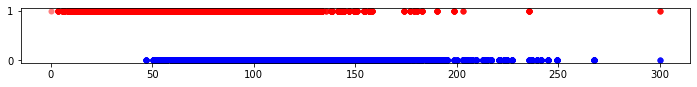}\\ \includegraphics[scale=0.5]{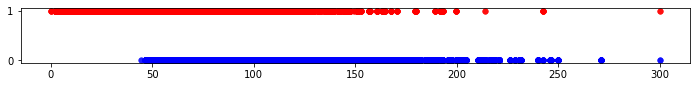}\\\includegraphics[scale=0.5]{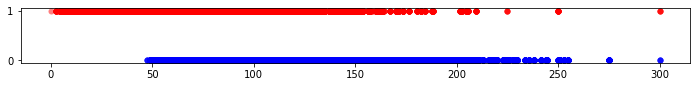}
\\\includegraphics[scale=0.5]{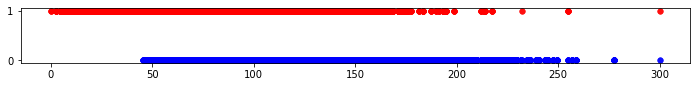}
\end{center}
\caption{\small\sf Ranges of $b$-loadings for $n=10,11,12,13$ from top to bottom. A red (resp. blue) dot at $(x,1)$ (resp. $(x,0)$) corresponds to $\mathbf t \in \mathcal P(n)^3$ with $b(\mathbf t) =x$ and $g(\mathbf t) \neq 0$ (resp. $g(\mathbf t)=0$).} \label{fig-2} \end{figure}

\begin{figure}[h!]\begin{center}
\includegraphics[scale=0.35]{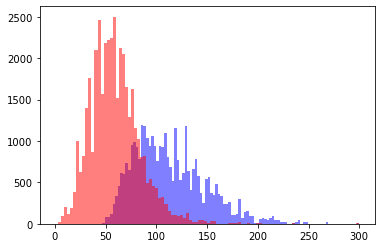}\quad \includegraphics[scale=0.35]{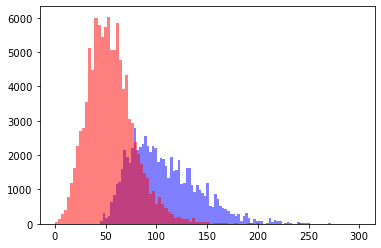}\\ \includegraphics[scale=0.35]{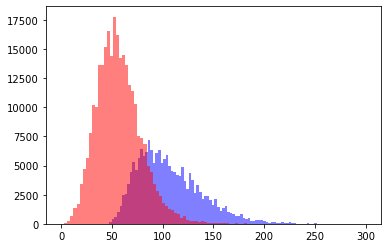}
\quad \includegraphics[scale=0.35]{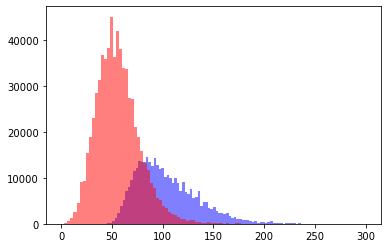}
\end{center}
\caption{\small \sf Histograms of $b$-loadings for $n=10$ (top-left), $11$ (top-right), $12$ (bottom-left) and $13$ (bottom-right). The red (resp. blue) region represents the numbers of $\mathbf t$ such that $g(\mathbf t) \neq 0$ (resp. $g(\mathbf t) =0$).} \label{fig-3} \end{figure}

In Figures \ref{fig-4}--\ref{fig-3}, we present the ranges and histograms of loadings of $\mathbf t \in \mathcal P(n)^3$ for $n=10,11,12,13$ according to whether $g(\mathbf t) \neq 0$ (red) or $= 0$ (blue). 
 In Figure \ref{fig-4}, the $y$-values $0$ and $1$ represent the cases $g(\mathbf t) =0$ and $g(\mathbf t) \neq 0$, respectively, while the $x$-value is the $a$-loading of $\mathbf t$. The same convention applies to Figure \ref{fig-2} with $b$-loadings.
As one can see, the ranges and histograms do not vary much as $n$ varies. The separation between the regions corresponding to $g(\mathbf t) \neq 0$ (red) and $= 0$ (blue) is more distinctive in the case of $b$-loadings. It is clear that we may use the minimum values of loadings to obtain vertical lines that separate the red regions from the blue ones.

With this in mind, define
\begin{align*} a_{\star} & \seteq \min \{ a(\mathbf t) : g(\mathbf t) \neq 0, \mathbf t \in \mathcal P(n)^3 \}, \\  b_{\star} & \seteq \min \{ b(\mathbf t) : g(\mathbf t) = 0 , \mathbf t \in \mathcal P(n)^3 \}  . \end{align*}
Then, for $\mathbf t \in \mathcal P(n)^3$,
\begin{equation} \label{eqn-box}   \text{ if } a(\mathbf t) < a_{\star} \text{ then } g(\mathbf t) = 0 \ \text{ and }\ \boxed{\text{if } b(\mathbf t) < b_{\star} \text{ then } g(\mathbf t) \neq 0} . \end{equation}
This provides sufficient conditions for $g(\mathbf t) = 0$  and $g(\mathbf t) \neq 0$, respectively, once we know the values of $a_\star$ and $b_\star$.

 In this way, the values of $b_\star$ can be used quite effectively in distinguishing $g(\mathbf t) \neq 0$ from $g(\mathbf t) =0$.  When $n=20$, the percentage of $\mathbf t$ satisfying $b(\mathbf t) < b_\star$ is about 31.8\%. In contrast, the values $a_\star$ do not turn out to be very useful for bigger $n$ in distinguishing $g(\mathbf t) = 0$ from $g(\mathbf t) \neq 0$. When $n=20$, the percentage of $\mathbf t$ satisfying $a(\mathbf t) < a_\star$ is only $0.37\%$. See Example \ref{exa-18}\,(2) below for more details. Nonetheless, the values of $a_\star$
  are interesting in their own right and can be valuable for analyzing the distribution of the $a$-loadings in relation to the Kronecker coefficients.  

\begin{example} \label{exa-18} \hfill
\begin{enumerate}
\item When $n=18$, we obtain $b_{\star} \approx 44.18$. Now that the $b$-loading of \[ \mathbf t=((12, 4, 2), (8, 4, 2, 2, 1, 1), (5, 4, 3, 3, 1, 1, 1))\] is readily computed to be approximately $41.07 < b_{\star}$, we immediately conclude that $g(\mathbf t) \neq 0$ by \eqref{eqn-box}.
\item When $n=20$, there are $246,491,883$ triples $\mathbf t \in \mathcal P(20)$. Among them, $78,382,890$ triples satisfy $b(\mathbf t) < b_\star \approx 43.74$ so that $g(\mathbf t) \neq 0$. The percentage of these triples is about 31.8\%. On the other hand, $909,200$ triples satisfy $a(\mathbf t) < a_\star \approx 70.88$ and the percentage is only $0.37\%$. 

\end{enumerate}
\end{example}

\begin{remark}
It appears that the $b$-loadings of $\mathbf t$ with $g(\mathbf t) \neq 0$ is a gamma distribution by itself. See the histogram and the curve of a gamma distribution when $n=13$ in Figure \ref{fig-br}.
\begin{figure}[h!]
\begin{center}
\includegraphics[scale=0.4]{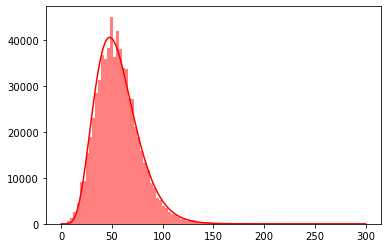}
\end{center}
\caption{\small\sf Histogram and curve (red) of a gamma distribution when $n=13$}\label{fig-br}
\end{figure}
\end{remark}

\medskip

In the rest of this section, computational results of the values of $a_{\star}$ and $b_{\star}$ for $6 \le n \le 20$ will be presented along with some conjectures. 

\subsection{Results on $a$-loadings}
We compute and record $a_{\star}$ and $\mathbf t=(\lambda, \mu, \nu) $ such that $a_{\star}=a(\mathbf t) $ and $\lambda \ge \mu \ge \nu$ lexicographically,  for $6 \le n \le 20$ in Table \ref{tab-2}. We do not consider $n \le 5$ because they seem to be too small for statistical analysis. 

\begin{table}[h!]
\begin{center}
\begin{tabular}{|c|c|l|l|l|} \hline
$n$ & $a_{\star}$  & \phantom{LL} $\lambda$ & \phantom{LLL} $\mu$ & \phantom{LLLL} $\nu$ \\ \hline
6 & 90.9986 & $(3,3)$&  $(2,2,2)$ & $(1,1,1,1,1,1)$  \\  \hline
7 & 85.0932 &$(2,2,2,1)$&$(2,2,2,1)$&$(2,2,2,1)$\\ \hline 
8 & 79.1637 & $(2^4)$  &  $(2^4)$  & $(2^4)$  \\ \hline
9 & 84.5605 &  $(3,2,2,2)$  &  $(2,2,2,2,1)$  &  $(2,2,2,2,1)$  \\ \hline
10 & 82.5959 & $(3,3,2,2)$& $(2,2,2,2,2)$&  $(2,2,2,2,1,1)$ \\ \hline
11& 78.1018 & $(3,3,3,2)$ & $(2^5,1)$ & $(2^5,1)$ \\ \hline
12& 74.6018 &  $(3^4)$   &  $(2^6)$  &  $(2^6)$  \\ \hline 
13& 78.1813& $(4,3,3,3)$& $(2^6,1)$& $(2^6,1)$ \\ \hline
14& 77.3651& $(4,4,3,3)$ &$(2^7)$&$(2^6,1,1)$ \\ \hline
15& 74.8437  &$(4,4,4,3)$&$(2^7,1)$&$(2^7,1)$  \\ \hline
16& 72.1837&  $(4^4)$  &  $(2^8)$  &  $(2^8)$  \\ \hline
17&71.2716  &$(3^5,2)$&$(3^5,2)$&$(2^8,1)$ \\ \hline
18& 68.9559 &$(3^6)$&$(3^6)$&$(2^9)$  \\ \hline
19& 71.9678 &$(4,3^5)$&$(3^6,1)$&$(2^9,1)$  \\ \hline
20& 70.8806 & $(5^4)$  &  $(2^{10})$  &  $(2^{10})$ \\ \hline
\end{tabular} 
\end{center}
\caption{\small\sf Values of $a_{\star}$ and $\mathbf t=(\lambda, \mu, \nu) $ such that $a_{\star}=a(\mathbf t) $ and $\lambda \ge \mu \ge \nu$ lexicographically.  When $n=8,12,16,20$, the partitions $\lambda, \mu, \nu$ are highlighted in blue to emphasize a pattern leading to Conjecture   \ref{conj-red}.  } \label{tab-2}
\end{table}
\medskip

Based on the results of $n=8,12,16,20$ as written in blue in Table \ref{tab-2}, we make the following conjecture.

\begin{conjecture} \label{conj-red}
Recall $a_{\star}  \seteq \min \{ a(\mathbf t) : g(\mathbf t) \neq 0, \mathbf t \in \mathcal P(n)^3 \}$, where $a(\mathbf t)\seteq a_\lambda+a_\mu+a_\nu$  and $g(\mathbf t) \seteq g_{\lambda, \mu}^\nu$ for $\mathbf t=(\lambda, \mu, \nu) \in \mathcal P(n)^3$. When $n=4k$ ($k \ge 2$), the values $a_{\star}$ are attained by $\mathbf t=((k^4), (2^{2k}), (2^{2k}))$.
\end{conjecture}

As an exhaustive computation for all possible triples becomes exponentially expensive, we assume that Conjecture \ref{conj-red} is true and continue computation. The results are in Table \ref{tab-3}.
Since we know $\mathbf t$ exactly under Conjecture \ref{conj-red}, we could calculate $a_{\star}$ for $n$ much bigger than those $n$ in the case of $b_{\star}$ that will be presented in Table \ref{tab-4}.

\begin{table}[h!] 
\begin{center}
\begin{tabular}{|c|c|l|} \hline
$n$ & $a_{\star}$  & \phantom{LLLLL} $\mathbf t$ \\ \hline
24& 70.0772 &$((6^4), (2^{12}), (2^{12}))$\\ \hline
28& 69.5351 &$((7^4), (2^{14}), (2^{14}))$ \\ \hline
32& 69.1732 &$((8^4), (2^{16}), (2^{16}))$\\ \hline 
36& 68.9254 &$((9^4), (2^{18}), (2^{18}))$\\ \hline 
40& 68.7518 &$((10^4), (2^{20}), (2^{20}))$\\ \hline 
44& 68.6334 &$((11^4), (2^{22}), (2^{22}))$\\ \hline 
48& 68.5549 &$((12^4), (2^{24}), (2^{24}))$\\ \hline 
\end{tabular} 
\end{center}
\caption{\small\sf Under Conjecture \ref{conj-red}, values of $a_{\star}$ and $\mathbf t=((k^4), (2^{2k}), (2^{2k}))$ for $n=4k$ such that $a_{\star}=a(\mathbf t) $ } \label{tab-3}
\end{table}

\begin{remark}
The values of $a_{\star}$ seem to keep  decreasing though slowly. However, it is not clear whether $a_\star$ converges to a limit as $n \to \infty$. 
\end{remark}

\medskip

\begin{figure}[h!] \begin{center}
\includegraphics[width=1.0\textwidth]{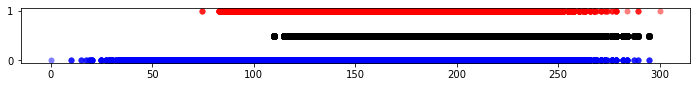}\\ \includegraphics[scale=0.45]{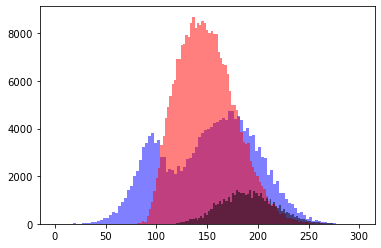}
\end{center} 
\caption{\small\sf Ranges of $a$-loadings where a black dot at $(x,\frac 1 2 )$ corresponds to $\mathbf t \in \mathcal P(12)^3$ satisfying the condition in \eqref{eqn-dl-1} with $a(\mathbf t) =x$ and $g(\mathbf t) = 0$ and histograms of $a$-loadings where the dark brown region  represents the numbers of $\mathbf t$ satisfying the condition in \eqref{eqn-dl-1} and $g(\mathbf t) = 0$.  } \label{fig-1} \end{figure}

Notice that we have a sufficient condition for $g(\mathbf t) =0$ by taking the contrapositive of \eqref{eqn-dl}:
\begin{equation} \label{eqn-dl-1}   d_\nu < |d_\lambda - d_\mu | \quad \text{ or } \quad  d_\nu > d_\lambda + d_\mu \quad \Longrightarrow  \quad  g(\mathbf t)= 0 . \end{equation}  
As $a_{\star}$ provides another sufficient condition for $g(\mathbf t) =0$ in \eqref{eqn-box}, one may be curious about their relationship. As a matter of fact, we observe that
\[  a_{\star} < a(\mathbf t) \quad \text{ for any $ \mathbf t$ satisfying the condition in \eqref{eqn-dl-1}}.\] Thus conditions in \eqref{eqn-box} and \eqref{eqn-dl-1} for $g(\mathbf t) =0$ do not have overlaps.
Let us look at pictures when $n=12$. In the  top  graph of Figure \ref{fig-1}, red dots and blue dots are   the same as in Figure \ref{fig-4}, while a black dot at $(x,\frac 1 2 )$ corresponds to $\mathbf t \in \mathcal P(12)^3$ satisfying the condition in \eqref{eqn-dl-1} with $a(\mathbf t) =x$ and $g(\mathbf t) = 0$.
In the  bottom  histograms of Figure \ref{fig-1}, the red region and blue region are  the same as in Figure \ref{fig-5}, while the dark brown region  represents the numbers of $\mathbf t$ satisfying the condition in \eqref{eqn-dl-1} and $g(\mathbf t) = 0$. 

\medskip

\subsection{Results on $b$-loadings}
 In Table \ref{tab-1}, we list $b_{\star}$ and $\mathbf t=(\lambda, \mu, \nu) $ such that $b_{\star}=b(\mathbf t) $ and $\lambda \ge \mu \ge \nu$ in lexicographic order for $6 \le n \le 20$.  When there are more than one $\mathbf t$ such that  $b_{\star}=b(\mathbf t)$, we only record the lexicographically smallest one. (Recall Remark \ref{quote-1}.)

\begin{example} 
When $n=16$, we get 
$b_{\star}=b(\mathbf t_1) = b(\mathbf t_2) = b(\mathbf t_3) $ with
\begin{align*} \mathbf t_1 &=[(10,3,2,1),(10,3,2,1),(5,3,2,1^6)],\\  \mathbf t_2&=[(10,3,2,1),(9,3,2,1,1),(4,3,2,1^7)], \\ \mathbf t_3& =[(5,3,2,1^6),(4,3,2,1^7),(4,3,2,1^7)], \end{align*} 
and only $\mathbf t_3$ is recorded in the table.
\end{example}

\begin{table}[h!]
\begin{center}
\begin{tabular}{|c|c|l|l|l|} \hline
$n$ & $b_{\star}$  & \phantom{LLLL} $\lambda$ & \phantom{LLLL} $\mu$ & \phantom{LLLL} $\nu$ \\ \hline
6 & 59.7812 &  $(2,2,1,1)$  &   $(2,2,1,1)$  &  $(2,2,1,1)$  \\  \hline
7 & 47.9477 &$(3,3,1)$&$(3,1,1,1,1)$&$(3,1,1,1,1)$\\ \hline 
8 & 54.6650 &$(3,2,1,1,1)$&$ (3,2,1,1,1)$&$(2,2,1,1,1,1)$\\ \hline
9 & 39.8213 &  $(3,2,1,1,1,1)$  &  $(3,2,1,1,1,1)$  &  $(3,2,1,1,1,1)$  \\ \hline
10 & 46.6592 & $(4,2,1^4)$& $(3,2,1^5)$&  $(3,2,1^5)$ \\ \hline
11& 44.4953 & $(6,1^5)$ & $(6,1^5)$ & $(4,3,3,1)$ \\ \hline
12& 47.3571 &  $(3,3,2,1^4)$  &  $(3,3,2,1^4)$  &  $(3,3,2,1^4)$  \\ \hline 
13& 45.1104& $(4,3,2, 1^4)$& $(3,3,2,1^5)$& $(3,3,2,1^5)$ \\ \hline
14& 44.9312& $(4,3,2,1^5)$ &$(4,3,2,1^5)$&$(3,3,2,1^6)$ \\ \hline
15& 40.3916  &  $(4,3,2,1^6)$  &  $(4,3,2,1^6)$  &  $(4,3,2,1^6)$  \\ \hline
16& 41.7064 & $(5,3,2,1^6)$ & $(4,3,2,1^7)$&$(4,3,2,1^7)$ \\ \hline
17& 43.4181 &$(5, 3, 2, 1^7)$&$(4,3,2,2,1^6)$&$(4,3,2,2,1^6)$\\ \hline
18& 44.1817 &   $(4,4,2,2,1^6)$  &   $(4,4,2,2,1^6)$ & $(4,4,2,2,1^6)$  \\ \hline
19& 44.3797 &$(5, 4, 2,2, 1^6)$&$(4,4,2,2,1^7)$&$(4,4,2,2,1^7)$\\ \hline
20& 43.7424 &$(5, 4, 2,2, 1^7)$&$(4,4,3,2,1^7)$&$(4,4,3,2,1^7)$\\ \hline
\end{tabular} 
\end{center}
\caption{\small\sf Values of $b_{\star}$ and $\mathbf t=(\lambda, \mu, \nu) $ such that $b_{\star}=b(\mathbf t) $.   When $n=6,9,12,15,18$, the partitions $\lambda, \mu, \nu$ are highlighted in blue to emphasize the pattern $\lambda = \mu =\nu$ leading to Conjecture  \ref{conj-blue}. }\label{tab-1}
\end{table}
%

 Drawing from the results in Table \ref{tab-1}---particularly those for $n=6,9,12,15,18$ highlighted in blue---we propose the following conjecture. 

\begin{conjecture} \label{conj-blue}
Recall $b_{\star}  \seteq \min \{ b(\mathbf t) : g(\mathbf t) = 0, \mathbf t \in \mathcal P(n)^3 \}$, where $b(\mathbf t)\seteq b_\lambda+b_\mu+b_\nu$  and $g(\mathbf t) \seteq g_{\lambda, \mu}^\nu$ for $\mathbf t=(\lambda, \mu, \nu) \in \mathcal P(n)^3$. For $n\ge 6$, the values $b_{\star}$ are attained by $\mathbf t=(\lambda,\mu,\nu)$ such that $\lambda=\mu$ or $\mu=\nu$. Moreover, when $n=3k$, $k \ge 2$, the values $b_{\star}$ are attained by $\mathbf t=(\lambda,\mu,\nu)$ such that $\lambda=\mu=\nu$.
\end{conjecture}

 Analogous to the case of $a_\star$, we assume Conjecture \ref{conj-blue} holds for  $n=3k$ and proceed with the computation. The results are presented in Table  \ref{tab-4}. 

\begin{table}[h!] 
\begin{center}
\begin{tabular}{|c|c|l|} \hline
$n$ & $b_{\star}$  & \phantom{LL} $\lambda=\mu=\nu$ \\ \hline
21& 45.0545 &$(5,4,2,2,1^8)$\\ \hline
24& 43.7126 &$(5,4,3,2,2,1^8)$ \\ \hline
27& 44.0699 &$(5,5,3,3,2,1^9)$\\ \hline 
30& 45.0141 &$(5,5,4,3,2,2,1^9)$\\ \hline 
33& 44.7615 &$(6,6,4,3,2,1^{12})$\\ \hline 
36& 44.3350 &$(6,6,4,3,2^3,1^{11})$\\ \hline 
\end{tabular} 
\end{center}
\caption{\small\sf Under Conjecture \ref{conj-blue}, values of $b_{\star}$ and $\mathbf t=(\lambda, \lambda, \lambda)$ for $n=3k$ such that $b_{\star}=b(\mathbf t)$} \label{tab-4}
\end{table}

\begin{remark}
The values of $b_{\star}$ appear to fluctuate with diminishing amplitudes as $n$ increases. However, it remains unclear whether $b_{\star}$ converges as $n \to \infty$. 
\end{remark}

\pagebreak

\appendix

\renewcommand{\baselinestretch}{1.0}
\section{Table of Loadings} \label{appen}
We tabulate the $a$-loading $a_\lambda$ and $b$-loading $b_\lambda$ of each partition $\lambda \in \mathcal P(n)$ for $6 \le n \le 12$.
{\scriptsize
\begin{center}
\begin{tabular}{|c|r|r|}\hline
$\lambda$ & $a_\lambda$\phantom{LL} & $b_\lambda$\phantom{LL}  \\ \hline
 $(6, 0, 0, 0, 0, 0)$& 100.0 & 100.0 \\
 $(5, 1, 0, 0, 0, 0)$& 85.8934& 37.252\\
 $(4, 2, 0, 0, 0, 0)$& 71.7868& 19.9271 \\
 $(4, 1, 1, 0, 0, 0)$& 66.6591& 4.363\\
 $(3, 3, 0, 0, 0, 0)$& 57.6803& 43.005\\
 $(3, 2, 1, 0, 0, 0)$& 52.5526& 0.0 \\
 $(3, 1, 1, 1, 0, 0)$& 45.2311& 4.363\\
 $(2, 2, 2, 0, 0, 0)$& 33.3183& 43.005\\
 $(2, 2, 1, 1, 0, 0)$& 31.1245& 19.9271 \\
 $(2, 1, 1, 1, 1, 0)$& 22.8133& 37.252\\
 $(1, 1, 1, 1, 1, 1)$& 0.0 & 100.0 \\ \hline
$(7, 0, 0, 0, 0, 0, 0)$& 100.0& 100.0 \\
 $(6, 1, 0, 0, 0, 0, 0)$& 88.302& 47.507\\
 $(5, 2, 0, 0, 0, 0, 0)$& 76.604& 26.483\\
 $(5, 1, 1, 0, 0, 0, 0)$& 72.8338& 13.1061\\
 $(4, 3, 0, 0, 0, 0, 0)$& 64.906& 36.928\\
 $(4, 2, 1, 0, 0, 0, 0)$& 61.1358& 0.0\\
 $(4, 1, 1, 1, 0, 0, 0)$& 55.5306& 1.81\\
 $(3, 3, 1, 0, 0, 0, 0)$& 49.4378& 21.735\\
 $(3, 2, 2, 0, 0, 0, 0)$& 45.6676& 21.735\\
 $(3, 2, 1, 1, 0, 0, 0)$& 43.8326& 0.0\\
 $(3, 1, 1, 1, 1, 0, 0)$& 37.3978& 13.1061\\
 $(2, 2, 2, 1, 0, 0, 0)$& 28.3644& 36.928\\
 $(2, 2, 1, 1, 1, 0, 0)$& 25.6998& 26.483\\
 $(2, 1, 1, 1, 1, 1, 0)$& 18.7933& 47.507\\
 $(1, 1, 1, 1, 1, 1, 1)$& 0.0& 100.0 \\ \hline
$(8, 0, 0, 0, 0, 0, 0, 0)$& 100.0& 100.0\\
 $(7, 1, 0, 0, 0, 0, 0, 0)$& 90.5921& 58.055\\
 $(6, 2, 0, 0, 0, 0, 0, 0)$& 81.1842& 35.198\\
 $(6, 1, 1, 0, 0, 0, 0, 0)$& 77.6539& 28.246\\
 $(5, 3, 0, 0, 0, 0, 0, 0)$& 71.7763& 34.854\\
 $(5, 2, 1, 0, 0, 0, 0, 0)$& 68.2461& 9.7331\\
 $(5, 1, 1, 1, 0, 0, 0, 0)$& 63.194& 12.637\\
 $(4, 4, 0, 0, 0, 0, 0, 0)$& 62.3685& 48.552\\
 $(4, 3, 1, 0, 0, 0, 0, 0)$& 58.8382& 15.265\\
 $(4, 2, 2, 0, 0, 0, 0, 0)$& 55.3079& 15.265\\
 $(4, 2, 1, 1, 0, 0, 0, 0)$& 53.7861& 0.0\\
 $(4, 1, 1, 1, 1, 0, 0, 0)$& 47.8449& 12.637\\
 $(3, 3, 2, 0, 0, 0, 0, 0)$& 45.9& 30.531\\
 $(3, 3, 1, 1, 0, 0, 0, 0)$& 44.3782& 15.265\\
 $(3, 2, 2, 1, 0, 0, 0, 0)$& 40.8479& 15.265\\
 $(3, 2, 1, 1, 1, 0, 0, 0)$& 38.437& 9.7331\\
 $(3, 1, 1, 1, 1, 1, 0, 0)$& 32.0837& 28.246\\
 $(2, 2, 2, 2, 0, 0, 0, 0)$& 26.3879& 48.552\\
 $(2, 2, 2, 1, 1, 0, 0, 0)$& 25.4988& 34.854\\
 $(2, 2, 1, 1, 1, 1, 0, 0)$& 22.6758& 35.198\\
 $(2, 1, 1, 1, 1, 1, 1, 0)$& 16.0886& 58.055\\
 $(1, 1, 1, 1, 1, 1, 1, 1)$& 0.0& 100.0\\ \hline
$(9, 0, 0, 0, 0, 0, 0, 0, 0)$& 100.0& 100.0\\
 $(8, 1, 0, 0, 0, 0, 0, 0, 0)$& 91.876& 62.802\\
 $(7, 2, 0, 0, 0, 0, 0, 0, 0)$& 83.7521& 39.559\\
 $(7, 1, 1, 0, 0, 0, 0, 0, 0)$& 80.9205& 33.587\\
 $(6, 3, 0, 0, 0, 0, 0, 0, 0)$& 75.6281& 34.591\\
 $(6, 2, 1, 0, 0, 0, 0, 0, 0)$& 72.7965& 13.273\\
 $(6, 1, 1, 1, 0, 0, 0, 0, 0)$& 68.4825& 16.425\\
 $(5, 4, 0, 0, 0, 0, 0, 0, 0)$& 67.5041& 42.455\\
 $(5, 3, 1, 0, 0, 0, 0, 0, 0)$& 64.6726& 12.1941\\
 $(5, 2, 2, 0, 0, 0, 0, 0, 0)$& 61.841& 12.1941\\
 \hline
\end{tabular}
\begin{tabular}{|c|r|r|}\hline
$\lambda$ & $a_\lambda$\phantom{LL} & $b_\lambda$\phantom{LL}\\ \hline
 $(5, 2, 1, 1, 0, 0, 0, 0, 0)$& 60.3586& 0.0\\
 $(5, 1, 1, 1, 1, 0, 0, 0, 0)$& 55.3152& 10.278\\
 $(4, 4, 1, 0, 0, 0, 0, 0, 0)$& 56.5486& 26.205\\
 $(4, 3, 2, 0, 0, 0, 0, 0, 0)$& 53.7171& 17.261\\
 $(4, 3, 1, 1, 0, 0, 0, 0, 0)$& 52.2346& 5.067\\
 $(4, 2, 2, 1, 0, 0, 0, 0, 0)$& 49.4031& 5.067\\
 $(4, 2, 1, 1, 1, 0, 0, 0, 0)$& 47.1912& 0.0\\
 $(4, 1, 1, 1, 1, 1, 0, 0, 0)$& 41.7289& 16.425\\
 $(3, 3, 3, 0, 0, 0, 0, 0, 0)$& 42.7616& 39.778\\
 $(3, 3, 2, 1, 0, 0, 0, 0, 0)$& 41.2791& 17.261\\
 $(3, 3, 1, 1, 1, 0, 0, 0, 0)$& 39.0672& 12.1941\\
 $(3, 2, 2, 2, 0, 0, 0, 0, 0)$& 36.9651& 26.205\\
 $(3, 2, 2, 1, 1, 0, 0, 0, 0)$& 36.2357& 12.1941\\
 $(3, 2, 1, 1, 1, 1, 0, 0, 0)$& 33.6049& 13.273\\
 $(3, 1, 1, 1, 1, 1, 1, 0, 0)$& 27.9202& 33.587\\
 $(2, 2, 2, 2, 1, 0, 0, 0, 0)$& 23.7977& 42.455\\
 $(2, 2, 2, 1, 1, 1, 0, 0, 0)$& 22.6494& 34.591\\
 $(2, 2, 1, 1, 1, 1, 1, 0, 0)$& 19.7962& 39.559\\
 $(2, 1, 1, 1, 1, 1, 1, 1, 0)$& 13.9854& 62.802\\
 $(1, 1, 1, 1, 1, 1, 1, 1, 1)$& 0.0& 100.0\\ \hline
$(10, 0, 0, 0, 0, 0, 0, 0, 0, 0)$& 100.0& 100.0\\
 $(9, 1, 0, 0, 0, 0, 0, 0, 0, 0)$& 93.0766& 67.7441\\
 $(8, 2, 0, 0, 0, 0, 0, 0, 0, 0)$& 86.1532& 45.12\\
 $(8, 1, 1, 0, 0, 0, 0, 0, 0, 0)$& 83.5036& 41.476\\
 $(7, 3, 0, 0, 0, 0, 0, 0, 0, 0)$& 79.2298& 36.947\\
 $(7, 2, 1, 0, 0, 0, 0, 0, 0, 0)$& 76.5802& 20.739\\
 $(7, 1, 1, 1, 0, 0, 0, 0, 0, 0)$& 72.6788& 23.437\\
 $(6, 4, 0, 0, 0, 0, 0, 0, 0, 0)$& 72.3065& 39.542\\
 $(6, 3, 1, 0, 0, 0, 0, 0, 0, 0)$& 69.6568& 15.044\\
 $(6, 2, 2, 0, 0, 0, 0, 0, 0, 0)$& 67.0072& 15.044\\
 $(6, 2, 1, 1, 0, 0, 0, 0, 0, 0)$& 65.7554& 5.179\\
 $(6, 1, 1, 1, 1, 0, 0, 0, 0, 0)$& 61.1395& 14.455\\
 $(5, 5, 0, 0, 0, 0, 0, 0, 0, 0)$& 65.3831& 49.1901\\
 $(5, 4, 1, 0, 0, 0, 0, 0, 0, 0)$& 62.7334& 21.441\\
 $(5, 3, 2, 0, 0, 0, 0, 0, 0, 0)$& 60.0838& 13.151\\
 $(5, 3, 1, 1, 0, 0, 0, 0, 0, 0)$& 58.832& 3.286\\
 $(5, 2, 2, 1, 0, 0, 0, 0, 0, 0)$& 56.1824& 3.286\\
 $(5, 2, 1, 1, 1, 0, 0, 0, 0, 0)$& 54.2161& 0.0\\
 $(5, 1, 1, 1, 1, 1, 0, 0, 0, 0)$& 49.2041& 14.455\\
 $(4, 4, 2, 0, 0, 0, 0, 0, 0, 0)$& 53.1604& 24.72\\
 $(4, 4, 1, 1, 0, 0, 0, 0, 0, 0)$& 51.9086& 14.862\\
 $(4, 3, 3, 0, 0, 0, 0, 0, 0, 0)$& 50.5108& 27.307\\
 $(4, 3, 2, 1, 0, 0, 0, 0, 0, 0)$& 49.259& 6.572\\
 $(4, 3, 1, 1, 1, 0, 0, 0, 0, 0)$& 47.2927& 3.286\\
 $(4, 2, 2, 2, 0, 0, 0, 0, 0, 0)$& 45.3575& 14.862\\
 $(4, 2, 2, 1, 1, 0, 0, 0, 0, 0)$& 44.6431& 3.286\\
 $(4, 2, 1, 1, 1, 1, 0, 0, 0, 0)$& 42.2807& 5.179\\
 $(4, 1, 1, 1, 1, 1, 1, 0, 0, 0)$& 37.0369& 23.437\\
 $(3, 3, 3, 1, 0, 0, 0, 0, 0, 0)$& 39.686& 27.307\\
 $(3, 3, 2, 2, 0, 0, 0, 0, 0, 0)$& 38.4341& 24.72\\
 $(3, 3, 2, 1, 1, 0, 0, 0, 0, 0)$& 37.7197& 13.151\\
 $(3, 3, 1, 1, 1, 1, 0, 0, 0, 0)$& 35.3573& 15.044\\
 $(3, 2, 2, 2, 1, 0, 0, 0, 0, 0)$& 33.8182& 21.441\\
 $(3, 2, 2, 1, 1, 1, 0, 0, 0, 0)$& 32.7077& 15.044\\
 $(3, 2, 1, 1, 1, 1, 1, 0, 0, 0)$& 30.1135& 20.739\\
 $(3, 1, 1, 1, 1, 1, 1, 1, 0, 0)$& 24.747& 41.476\\
 $(2, 2, 2, 2, 2, 0, 0, 0, 0, 0)$& 22.2789& 49.1901\\
 $(2, 2, 2, 2, 1, 1, 0, 0, 0, 0)$& 21.8828& 39.542\\ \hline
\end{tabular}
\end{center}
\begin{center}
\begin{tabular}{|c|r|r|}\hline
$\lambda$ & $a_\lambda$\phantom{LL} & $b_\lambda$\phantom{LL} \\ \hline
 $(2, 2, 2, 1, 1, 1, 1, 0, 0, 0)$& 20.5405& 36.947\\
 $(2, 2, 1, 1, 1, 1, 1, 1, 0, 0)$& 17.8237& 45.12\\
 $(2, 1, 1, 1, 1, 1, 1, 1, 1, 0)$& 12.3875& 67.7441\\
 $(1, 1, 1, 1, 1, 1, 1, 1, 1, 1)$& 0.0& 100.0\\ \hline
$(11, 0, 0, 0, 0, 0, 0, 0, 0, 0, 0)$& 100.0& 100.0\\
 $(10, 1, 0, 0, 0, 0, 0, 0, 0, 0, 0)$& 93.8295& 71.265\\
 $(9, 2, 0, 0, 0, 0, 0, 0, 0, 0, 0)$& 87.6591& 49.697\\
 $(9, 1, 1, 0, 0, 0, 0, 0, 0, 0, 0)$& 85.397& 46.624\\
 $(8, 3, 0, 0, 0, 0, 0, 0, 0, 0, 0)$& 81.4886& 39.924\\
 $(8, 2, 1, 0, 0, 0, 0, 0, 0, 0, 0)$& 79.2265& 26.3731\\
 $(8, 1, 1, 1, 0, 0, 0, 0, 0, 0, 0)$& 75.8034& 28.711\\
 $(7, 4, 0, 0, 0, 0, 0, 0, 0, 0, 0)$& 75.3182& 39.780\\
 $(7, 3, 1, 0, 0, 0, 0, 0, 0, 0, 0)$& 73.0561& 18.329\\
 $(7, 2, 2, 0, 0, 0, 0, 0, 0, 0, 0)$& 70.794& 18.329\\
 $(7, 2, 1, 1, 0, 0, 0, 0, 0, 0, 0)$& 69.6329& 10.1901\\
 $(7, 1, 1, 1, 1, 0, 0, 0, 0, 0, 0)$& 65.58& 17.872\\
 $(6, 5, 0, 0, 0, 0, 0, 0, 0, 0, 0)$& 69.1477& 45.913\\
 $(6, 4, 1, 0, 0, 0, 0, 0, 0, 0, 0)$& 66.8856& 20.801\\
 $(6, 3, 2, 0, 0, 0, 0, 0, 0, 0, 0)$& 64.6236& 12.90\\
 $(6, 3, 1, 1, 0, 0, 0, 0, 0, 0, 0)$& 63.4625& 4.762\\
 $(6, 2, 2, 1, 0, 0, 0, 0, 0, 0, 0)$& 61.2004& 4.762\\
 $(6, 2, 1, 1, 1, 0, 0, 0, 0, 0, 0)$& 59.4095& 1.967\\
 $(6, 1, 1, 1, 1, 1, 0, 0, 0, 0, 0)$& 54.969& 14.412\\
 $(5, 5, 1, 0, 0, 0, 0, 0, 0, 0, 0)$& 60.7152& 30.394\\
 $(5, 4, 2, 0, 0, 0, 0, 0, 0, 0, 0)$& 58.4531& 18.834\\
 $(5, 4, 1, 1, 0, 0, 0, 0, 0, 0, 0)$& 57.292& 10.694\\
 $(5, 3, 3, 0, 0, 0, 0, 0, 0, 0, 0)$& 56.191& 21.014\\
 $(5, 3, 2, 1, 0, 0, 0, 0, 0, 0, 0)$& 55.0299& 2.795\\
 $(5, 3, 1, 1, 1, 0, 0, 0, 0, 0, 0)$& 53.2391& 0.0\\
 $(5, 2, 2, 2, 0, 0, 0, 0, 0, 0, 0)$& 51.6068& 10.694\\
 $(5, 2, 2, 1, 1, 0, 0, 0, 0, 0, 0)$& 50.977& 0.0\\
 $(5, 2, 1, 1, 1, 1, 0, 0, 0, 0, 0)$& 48.7986& 1.967\\
 $(5, 1, 1, 1, 1, 1, 1, 0, 0, 0, 0)$& 44.1465& 17.872\\
 $(4, 4, 3, 0, 0, 0, 0, 0, 0, 0, 0)$& 50.0206& 31.70\\
 $(4, 4, 2, 1, 0, 0, 0, 0, 0, 0, 0)$& 48.8595& 13.4\\
 $(4, 4, 1, 1, 1, 0, 0, 0, 0, 0, 0)$& 47.0686& 10.694\\
 $(4, 3, 3, 1, 0, 0, 0, 0, 0, 0, 0)$& 46.5974& 15.6\\
 $(4, 3, 2, 2, 0, 0, 0, 0, 0, 0, 0)$& 45.4363& 13.4\\
 $(4, 3, 2, 1, 1, 0, 0, 0, 0, 0, 0)$& 44.8065& 2.795\\
 $(4, 3, 1, 1, 1, 1, 0, 0, 0, 0, 0)$& 42.6281& 4.762\\
 $(4, 2, 2, 2, 1, 0, 0, 0, 0, 0, 0)$& 41.3834& 10.694\\
 $(4, 2, 2, 1, 1, 1, 0, 0, 0, 0, 0)$& 40.366& 4.762\\
 $(4, 2, 1, 1, 1, 1, 1, 0, 0, 0, 0)$& 37.9761& 10.1901\\
 $(4, 1, 1, 1, 1, 1, 1, 1, 0, 0, 0)$& 33.1883& 28.711\\
 $(3, 3, 3, 2, 0, 0, 0, 0, 0, 0, 0)$& 37.0038& 31.70\\
 $(3, 3, 3, 1, 1, 0, 0, 0, 0, 0, 0)$& 36.374& 21.014\\
 $(3, 3, 2, 2, 1, 0, 0, 0, 0, 0, 0)$& 35.2129& 18.834\\
 $(3, 3, 2, 1, 1, 1, 0, 0, 0, 0, 0)$& 34.1956& 12.90\\
 $(3, 3, 1, 1, 1, 1, 1, 0, 0, 0, 0)$& 31.8056& 18.329\\
 $(3, 2, 2, 2, 2, 0, 0, 0, 0, 0, 0)$& 31.1599& 30.394\\
 $(3, 2, 2, 2, 1, 1, 0, 0, 0, 0, 0)$& 30.7724& 20.801\\
 $(3, 2, 2, 1, 1, 1, 1, 0, 0, 0, 0)$& 29.5435& 18.329\\
 $(3, 2, 1, 1, 1, 1, 1, 1, 0, 0, 0)$& 27.0179& 26.3731\\
 $(3, 1, 1, 1, 1, 1, 1, 1, 1, 0, 0)$& 22.1582& 46.624\\
 $(2, 2, 2, 2, 2, 1, 0, 0, 0, 0, 0)$& 20.549& 45.913\\
 $(2, 2, 2, 2, 1, 1, 1, 0, 0, 0, 0)$& 19.9499& 39.780\\
 $(2, 2, 2, 1, 1, 1, 1, 1, 0, 0, 0)$& 18.5853& 39.924\\
 $(2, 2, 1, 1, 1, 1, 1, 1, 1, 0, 0)$& 15.9878& 49.697\\
 $(2, 1, 1, 1, 1, 1, 1, 1, 1, 1, 0)$& 11.0873& 71.265\\
 $(1, 1, 1, 1, 1, 1, 1, 1, 1, 1, 1)$& 0.0& 100.0 \\ \hline
$(12, 0, 0, 0, 0, 0, 0, 0, 0, 0, 0, 0)$& 100.0& 100.0\\
 $(11, 1, 0, 0, 0, 0, 0, 0, 0, 0, 0, 0)$& 94.5958& 74.832\\
 $(10, 2, 0, 0, 0, 0, 0, 0, 0, 0, 0, 0)$& 89.1916& 54.707\\
\hline \end{tabular}
\begin{tabular}{|c|r|r|}\hline
$\lambda$ & $a_\lambda$\phantom{LL} & $b_\lambda$\phantom{LL} \\ \hline
 $(10, 1, 1, 0, 0, 0, 0, 0, 0, 0, 0, 0)$& 87.0838& 52.743\\
 $(9, 3, 0, 0, 0, 0, 0, 0, 0, 0, 0, 0)$& 83.7874& 43.844\\
 $(9, 2, 1, 0, 0, 0, 0, 0, 0, 0, 0, 0)$& 81.6796& 33.490\\
 $(9, 1, 1, 1, 0, 0, 0, 0, 0, 0, 0, 0)$& 78.5079& 35.703\\
 $(8, 4, 0, 0, 0, 0, 0, 0, 0, 0, 0, 0)$& 78.3831& 41.257\\
 $(8, 3, 1, 0, 0, 0, 0, 0, 0, 0, 0, 0)$& 76.2754& 23.775\\
 $(8, 2, 2, 0, 0, 0, 0, 0, 0, 0, 0, 0)$& 74.1676& 23.775\\
 $(8, 2, 1, 1, 0, 0, 0, 0, 0, 0, 0, 0)$& 73.1037& 17.598\\
 $(8, 1, 1, 1, 1, 0, 0, 0, 0, 0, 0, 0)$& 69.3148& 24.48\\
 $(7, 5, 0, 0, 0, 0, 0, 0, 0, 0, 0, 0)$& 72.9789& 44.246\\
 $(7, 4, 1, 0, 0, 0, 0, 0, 0, 0, 0, 0)$& 70.8711& 22.913\\
 $(7, 3, 2, 0, 0, 0, 0, 0, 0, 0, 0, 0)$& 68.7634& 15.785\\
 $(7, 3, 1, 1, 0, 0, 0, 0, 0, 0, 0, 0)$& 67.6995& 9.608\\
 $(7, 2, 2, 1, 0, 0, 0, 0, 0, 0, 0, 0)$& 65.5917& 9.608\\
 $(7, 2, 1, 1, 1, 0, 0, 0, 0, 0, 0, 0)$& 63.9106& 8.104\\
 $(7, 1, 1, 1, 1, 1, 0, 0, 0, 0, 0, 0)$& 59.7695& 18.845\\
 $(6, 6, 0, 0, 0, 0, 0, 0, 0, 0, 0, 0)$& 67.5747& 50.894\\
 $(6, 5, 1, 0, 0, 0, 0, 0, 0, 0, 0, 0)$& 65.4669& 28.138\\
 $(6, 4, 2, 0, 0, 0, 0, 0, 0, 0, 0, 0)$& 63.3591& 17.15\\
 $(6, 4, 1, 1, 0, 0, 0, 0, 0, 0, 0, 0)$& 62.2953& 10.981\\
 $(6, 3, 3, 0, 0, 0, 0, 0, 0, 0, 0, 0)$& 61.2514& 19.530\\
 $(6, 3, 2, 1, 0, 0, 0, 0, 0, 0, 0, 0)$& 60.1875& 3.854\\
 $(6, 3, 1, 1, 1, 0, 0, 0, 0, 0, 0, 0)$& 58.5064& 2.350\\
 $(6, 2, 2, 2, 0, 0, 0, 0, 0, 0, 0, 0)$& 57.0159& 10.981\\
 $(6, 2, 2, 1, 1, 0, 0, 0, 0, 0, 0, 0)$& 56.3986& 2.350\\
 $(6, 2, 1, 1, 1, 1, 0, 0, 0, 0, 0, 0)$& 54.3653& 4.700\\
 $(6, 1, 1, 1, 1, 1, 1, 0, 0, 0, 0, 0)$& 49.9969& 18.845\\
 $(5, 5, 2, 0, 0, 0, 0, 0, 0, 0, 0, 0)$& 57.9549& 25.7881\\
 $(5, 5, 1, 1, 0, 0, 0, 0, 0, 0, 0, 0)$& 56.8911& 19.6111\\
 $(5, 4, 3, 0, 0, 0, 0, 0, 0, 0, 0, 0)$& 55.8471& 24.3081\\
 $(5, 4, 2, 1, 0, 0, 0, 0, 0, 0, 0, 0)$& 54.7833& 8.631\\
 $(5, 4, 1, 1, 1, 0, 0, 0, 0, 0, 0, 0)$& 53.1022& 7.127\\
 $(5, 3, 3, 1, 0, 0, 0, 0, 0, 0, 0, 0)$& 52.6755& 11.003\\
 $(5, 3, 2, 2, 0, 0, 0, 0, 0, 0, 0, 0)$& 51.6117& 8.631\\
 $(5, 3, 2, 1, 1, 0, 0, 0, 0, 0, 0, 0)$& 50.9944& 0.0\\
 $(5, 3, 1, 1, 1, 1, 0, 0, 0, 0, 0, 0)$& 48.9611& 2.350\\
 $(5, 2, 2, 2, 1, 0, 0, 0, 0, 0, 0, 0)$& 47.8227& 7.127\\
 $(5, 2, 2, 1, 1, 1, 0, 0, 0, 0, 0, 0)$& 46.8533& 2.350\\
 $(5, 2, 1, 1, 1, 1, 1, 0, 0, 0, 0, 0)$& 44.5927& 8.104\\
 $(5, 1, 1, 1, 1, 1, 1, 1, 0, 0, 0, 0)$& 40.099& 24.48\\
 $(4, 4, 4, 0, 0, 0, 0, 0, 0, 0, 0, 0)$& 48.3351& 38.25\\
 $(4, 4, 3, 1, 0, 0, 0, 0, 0, 0, 0, 0)$& 47.2713& 19.634\\
 $(4, 4, 2, 2, 0, 0, 0, 0, 0, 0, 0, 0)$& 46.2074& 17.2631\\
 $(4, 4, 2, 1, 1, 0, 0, 0, 0, 0, 0, 0)$& 45.5902& 8.631\\
 $(4, 4, 1, 1, 1, 1, 0, 0, 0, 0, 0, 0)$& 43.5569& 10.981\\
 $(4, 3, 3, 2, 0, 0, 0, 0, 0, 0, 0, 0)$& 44.0997& 19.634\\
 $(4, 3, 3, 1, 1, 0, 0, 0, 0, 0, 0, 0)$& 43.4824& 11.003\\
 $(4, 3, 2, 2, 1, 0, 0, 0, 0, 0, 0, 0)$& 42.4185& 8.631\\
 $(4, 3, 2, 1, 1, 1, 0, 0, 0, 0, 0, 0)$& 41.4491& 3.854\\
 $(4, 3, 1, 1, 1, 1, 1, 0, 0, 0, 0, 0)$& 39.1885& 9.608\\
 $(4, 2, 2, 2, 2, 0, 0, 0, 0, 0, 0, 0)$& 38.6296& 19.6111\\
 $(4, 2, 2, 2, 1, 1, 0, 0, 0, 0, 0, 0)$& 38.2774& 10.981\\
 $(4, 2, 2, 1, 1, 1, 1, 0, 0, 0, 0, 0)$& 37.0807& 9.608\\
 $(4, 2, 1, 1, 1, 1, 1, 1, 0, 0, 0, 0)$& 34.6948& 17.598\\
 $(4, 1, 1, 1, 1, 1, 1, 1, 1, 0, 0, 0)$& 30.1207& 35.703\\
 $(3, 3, 3, 3, 0, 0, 0, 0, 0, 0, 0, 0)$& 35.5238& 38.25\\
 $(3, 3, 3, 2, 1, 0, 0, 0, 0, 0, 0, 0)$& 34.9065& 24.3081\\
 $(3, 3, 3, 1, 1, 1, 0, 0, 0, 0, 0, 0)$& 33.9371& 19.530\\
 $(3, 3, 2, 2, 2, 0, 0, 0, 0, 0, 0, 0)$& 33.2254& 25.7881\\
 $(3, 3, 2, 2, 1, 1, 0, 0, 0, 0, 0, 0)$& 32.8732& 17.15\\
 $(3, 3, 2, 1, 1, 1, 1, 0, 0, 0, 0, 0)$& 31.6765& 15.785\\
 $(3, 3, 1, 1, 1, 1, 1, 1, 0, 0, 0, 0)$& 29.2906& 23.775\\
 $(3, 2, 2, 2, 2, 1, 0, 0, 0, 0, 0, 0)$& 29.0843& 28.138\\
 \hline
\end{tabular}
\end{center}

\hskip 1.4 cm \begin{tabular}{|c|r|r|}\hline
$\lambda$ & $a_\lambda$\phantom{LL} & $b_\lambda$\phantom{LL} \\ \hline
 $(3, 2, 2, 2, 1, 1, 1, 0, 0, 0, 0, 0)$& 28.5049& 22.913\\
 $(3, 2, 2, 1, 1, 1, 1, 1, 0, 0, 0, 0)$& 27.1828& 23.775\\
 $(3, 2, 1, 1, 1, 1, 1, 1, 1, 0, 0, 0)$& 24.7165& 33.490\\
 $(3, 1, 1, 1, 1, 1, 1, 1, 1, 1, 0, 0)$& 20.0998& 52.743\\
 $(2, 2, 2, 2, 2, 2, 0, 0, 0, 0, 0, 0)$& 19.539& 50.894\\
 $(2, 2, 2, 2, 2, 1, 1, 0, 0, 0, 0, 0)$& 19.3117& 44.246\\
 $(2, 2, 2, 2, 1, 1, 1, 1, 0, 0, 0, 0)$& 18.6069& 41.257\\
 $(2, 2, 2, 1, 1, 1, 1, 1, 1, 0, 0, 0)$& 17.2045& 43.844\\
 $(2, 2, 1, 1, 1, 1, 1, 1, 1, 1, 0, 0)$& 14.6956& 54.707\\
 $(2, 1, 1, 1, 1, 1, 1, 1, 1, 1, 1, 0)$& 10.0548& 74.832\\
 $(1, 1, 1, 1, 1, 1, 1, 1, 1, 1, 1, 1)$& 0.0& 100.0 \\ \hline
\end{tabular} 
}

\end{document}